\documentstyle{amsart}

\theoremstyle{plain}

\newtheorem{thm}{Theorem}[section]
\newtheorem{lem}{Lemma}[section]
\newtheorem{cor}{Corollary}[section]
\newtheorem{prop}{Proposition}[section]

\theoremstyle{definition}

\theoremstyle{remark}

\newtheorem{rem}{Remark}[section] 
\newtheorem{ack}{Acknowledgments}  
\newtheorem{step}{Step}  

\numberwithin{equation}{section}

\newsymbol\twoheadrightarrow 1310

\newcommand{\twr}{\twoheadrightarrow}
\newcommand{\riso}{\stackrel{\sim}{\longrightarrow}}
\newcommand{\rank}{\operatorname{rank}}
\newcommand{\codim}{\operatorname{codim}}
\newcommand{\pic}{\operatorname{Pic}}

\newcommand{\Hom}{\operatorname{Hom}}
\newcommand{\ext}{\operatorname{Ext}}
\newcommand{\morph}{\operatorname{om}}
\newcommand{\col}{\colon}
\newcommand{\sm}{\setminus}
\newcommand{\im}{{\imath}}
\newcommand{\jm}{{\jmath}}
\newcommand{\vp}{\varphi}
\newcommand{\D}{\Delta}
\newcommand{\G}{\Gamma}
\newcommand{\lr}{\longrightarrow}
\newcommand{\lf}{\longleftarrow}

\newcommand{\cj}{{\cal J}}
\newcommand{\ce}{{\cal E}}
\newcommand{\co}{{\cal O}}
\newcommand{\cs}{{\cal S}}
\newcommand{\cg}{{\cal G}}
\newcommand{\cf}{{\cal F}}
\newcommand{\cu}{{\cal U}}
\newcommand{\cv}{{\cal V}}

\newcommand{\cl}{{\cal L}}
\newcommand{\ct}{{\cal T}}
\newcommand{\cw}{{\cal W}}
\newcommand{\ch}{{\cal H}}
\newcommand{\bbC}{{\Bbb C}}

\newcommand{\bbN}{{\Bbb N}}
\newcommand{\bbP}{{\Bbb P}}
\newcommand{\bbQ}{{\Bbb Q}}
\newcommand{\bbZ}{{\Bbb Z}}
\newcommand{\bbR}{{\Bbb R}}
\newcommand{\bbH}{{\Bbb H}}
\newcommand{\balp}{\boldsymbol{\alpha}}

\newcommand{\comp}{{\scriptstyle{\circ}}}

\title{Intermediate Jacobians of Moduli Spaces} 
\author{Donu Arapura}
\address{Department of Mathematics\\
Purdue University\\
West Lafayette, IN 47907\\
U.S.A.}
\thanks{First author partially supported by the NSF}
\email{dvb@@math.purdue.edu}
\author{Pramathanath Sastry}
\address{The Mehta Research Institute \\
Chhatnag, Jhusi, Allahabad District \\
U.P., INDIA 221 506}
\email{pramath@@mri.ernet.in}
\date{\today}
\begin{document}
\maketitle

\section{Introduction}\label{s:intro}
We work throughout over the complex numbers $\bbC$, i.e. all schemes are over
$\bbC$ and all maps of schemes are maps of $\bbC$-schemes. A curve, unless
otherwise stated, is a smooth complete curve. Points mean geometric points.
We will, as is usual in such situations, toggle between the algebraic
and analytic categories without warning. 
For a quasi-projective algebraic variety $Y$, the (mixed) Hodge structure
associated with its $i$-th cohomology will be denoted $H^i(Y)$.

For a curve $X$, $\cs\cu_X(n,\,L)$ will denote the moduli space of
{\it semi-stable} vector bundles of rank $n$ and determinant $L$. The
smooth open subvariety defining the {\it stable locus} will be denoted
$\cs\cu_X^s(n,\,L)$. We assume familiarity with the basic facts about
such a moduli space as laid out, for example in 
\cite{css-drez},\,pp.\,51--52,\,VI.A (see also
Theorems 10, 17 and 18 of {\it loc.cit.}).  Our principal
result is the following  theorem\,:

\begin{thm}\label{thm:main} Let $X$ be a curve of genus $g\ge 3$, 
$n\ge 2$ an integer, and 
$L$ a  line bundle of degree $d$ on $X$ with $d$ odd if $g=3$ and 
$n=2$. Let $\cs^s=\cs\cu_X^s(n,\,L)$.
Then  
$H^3(\cs^s)$ is a pure Hodge structure of type
$\{(1,2),\,(2,1)\}$, and it carries a 
natural polarization making the intermediate Jacobian
$$J^2(\cs^s) = \frac{H^3(\cs^s,\,\bbC)}
{F^2+H^3(\cs^s,\,\bbZ)} $$
into a principally polarized abelian variety. There is an
isomorphism of principally polarized abelian varieties $J(X)\simeq 
J^2(\cs^s)$.
 \end{thm}

The word ``natural'' above has the following meaning: an isomorphism 
between any two $\cs^s$'s as above will induce an isomorphism on third
cohomology  which will respect the indicated polarizations.
As an immediate corollary, we obtain the following  Torelli theorem: 

\begin{cor}\label{cor:main} Let $X$ and $X'$ be curves of genus $g\ge 3$,
$L$ and $L'$ line bundles of degree $d$ on $X$ and $X'$ respectively, and
$n\ge 2$ an integer. If 
\begin{equation}\label{eqn:stable}
\cs\cu_X^s(n,\,L) \simeq \cs\cu_{X'}^s(n,\,L')
\end{equation}

or if 
\begin{equation}\label{eqn:semistable}
\cs\cu_X(n,\,L) \simeq \cs\cu_{X'}(n,\,L')
\end{equation}

then
$$
X \simeq X' ,
$$
except when $g=3, n=2, (n,\,d)\ne 1$.
\end{cor}

\begin{pf} Since $\cs\cu_X^s(n,\,L)$ (resp. $\cs\cu_{X'}^s(n,\,L')$)
is the smooth locus of $\cs\cu_X(n,\,L)$ (resp. $\cs\cu_{X'}(n,\,L')$),
therefore it is enough to assume \eqref{eqn:stable} holds.
By assumption $J^2(\cs\cu^s_X(n,\,L))\simeq J^2(\cs\cu^s_{X'}(n,\,L'))$ 
as polarized abelian varieties.
Therefore $J(X)\simeq J(X')$,
and the corollary follows from the usual Torelli theorem. 
\end{pf}

The theorem is new for $(n,\,d) \ne 1$ 
(the so called ``non-coprime case"). When $(n,\,d)=1$ (the ``coprime case"),
the theorem (and its corollary) 
has been proven by Narasimhan and Ramanan \cite{N-R}, Tyurin
\cite{Ty} and (for $n=2$) by Mumford and Newstead \cite{M-N}. In the
non-coprime case, Kouvidakis and Pantev \cite{KP} have proved the
above corollary under the assumption \eqref{eqn:semistable},
and in fact the full result can be deduced from this case. \footnote{In 
fact, the
the exceptional case in the corollary can be eliminated using
the results in \cite{KP}} However
the present line of reasoning is extremely natural,
and is of a rather different character from that of Kouvidakis and Pantev. In 
particular, 
Theorem\, \ref{thm:main} will not follow from their techniques.
In the special case where $n=2$ and $L= \co_X$, Balaji \cite{balaji}
has shown a similar Torelli type theorem
 for Seshadri's canonical desingularization 
$N \to \cs\cu_X(2, \co_X)$ \cite{css-desing} in the range $g >3$. 
\footnote{Balaji states the result for $g \ge 3$, but his proof seems to
work only for $g > 3$. (See Remark\,\ref{rmk:balaji}).} 
In the coprime case, the proofs in \cite{M-N}
and \cite{N-R} rely on the fact that $\cs\cu_X^s(n,\,L) = \cs\cu_X(n,\,L)$,
and hence $\cs\cu_X^s(n,\,L)$ is smooth projective, and most importantly
the product $X \times \cs\cu_X(n,\,L)$ possesses a Poincar{\'e} bundle. In the
non-coprime case $\cs\cu_X^s(n,\,L)$ is  not complete and a result of
Ramanan (see \cite{R}) says that there is no Poincar{\'e} bundle on
$X \times U$ for any Zariski open subset $U$ of $\cs\cu_X(n,\,L)$.

We concentrate primarily on the non-coprime case---the only remaining case of 
interest. Our strategy is to use a Hecke correspondence to relate the Hodge 
structure on $H^3(\cs\cu^s_X(n,\,L))$ to that on $H^1(X)$. To this extent 
our proof resembles Balaji's in \cite{balaji}. We are able to deduce more
than Balaji does by imposing a polarization (which varies well with
$\cs\cu_X^s(n,\,L)$) on the Hodge structure of $H^3(\cs\cu_X^s(n,\,L))$.
This construction of the polarization needs a version of Lefschetz's
Hyperplane Theorem (for quasi projective varieties. See 
Theorem\,\ref{thm:lefschetz}). There is however another approach
to the problem of polarization, which uses M. Saito's theory of 
polarizations on Hodge modules (see Remark\,\ref{rmk:saito}).
  
\section{The Main Ideas}\label{s:main-idea}
For the rest of the paper, we fix a curve $X$ of genus $g$, $n \in {\bbN}$,
$d \in {\bbZ}$  and a line bundle $L$ of degree
$d$ on $X$. Assume, as in the main theorem, that if $n=2$, then $g \ge 4$,
and that $g \ge 3$ otherwise.
We shall   assume, with one brief exception in step 3 below, 
that  $(n,\,d)\ne 1$.

We will also assume, for the rest of the paper, that $0 < d \le n$. This
involves no loss of generality, for $\cs\cu_X(n,\,L)$ is canonically
isomorphic to $\cs\cu_X(n,\,L\otimes\xi^n)$ for every line bundle
$\xi$ on $X$. Let $\cs = \cs\cu_X(n,\,L)$ and $\cs^s =
\cs\cu^s_X(n,\,L)$ and let $U\subseteq \cs$  be a smooth
open set containing $\cs^s$.

The broad strategy of our proof is as follows\,: Fix a set $\chi =
\{x^1,\ldots,\,x^{d-1}\} \subset X$ of $d-1$ distinct points.

\begin{step} First show that there are isomorphisms
(modulo torsion), depending only
on $(X,\,L,\,\chi)$, of Hodge structures

$$
\psi_{X,L,\chi}\col H^1(X)(-1) \riso
H^3(\cs^s)
$$
where $(-1)$ is the Tate twist. The isomorphism should vary well with the
data $(X,\,L,\,\chi)$. More precisely, 
suppose $\widetilde{X}\overset{h}{\to} T$ is a family of curves
of genus $g$, $\cl$ a line bundle on $\widetilde{X}$, whose restrictions
to the fibres of $h$ are of degree $d$, and $\widetilde{\chi}$ 
a set of $d-1$ mutually disjoint $T$-valued points on $\widetilde{X}$.
Let the specialization of 
$(\widetilde{X},\,\cl,\,\widetilde{\chi})$ at $t\in T$ be
$(X_t,\,L_t,\,\chi_t)$. 
Let $\widetilde{\cs^s}\overset{g}{\to} T$ be the resulting family 
$\{\cs\cu_{X_t}^s(n,\,L_t)\}$. Then there is an isomorphism
(modulo torsion) of variation of
Hodge structures
$$
\widetilde{\psi}\col R^1h_*\bbZ(-1) \riso R^3g_*\bbZ,
$$
which specializes at each $t\in T$ to $\psi_{X_t,L_t,\chi_t}$.
Note that $\psi_{X,L,\chi}$ gives an isomorphism of complex tori
$$
\vp_{X,L,\chi}\col J(X) \riso J^2(\cs^s)
$$
also varies well with $(X,\,L,\,\chi)$.
\end{step}
\begin{step} Find a (possibly nonprincipal) polarization
 $\Theta(\cs^s)$ on $J^2(\cs^s)$ which
depends only on $\cs^s$, and varies well with $\cs^s$. Let 
$\mu = \mu_{X,L,\chi}$ be the polarization on $J(X)$ induced by
$\Theta(\cs^s)$ and $\vp_{X,L,\chi}$.
\end{step}

\begin{step}
In this step we relax the above  assumptions, and no
longer insist that $(n,d) \ne 1$.
Suppose Steps 1 and 2 have been taken (see \cite{N-R} for the coprime
case).
 Theorem\,\ref{thm:main} will follow
by showing that there exists an integer $m$ such that 
$\frac{1}{m}\Theta$ 
is principal, and that $J^2(\cs^s)$ equipped with this polarization
is isomorphic to $J(X)$ with its canonical polarization. 
The essence of the argument will  be 
to show that any
natural polarization on $J(X)$ must be a multiple of the standard one.
  The argument is lifted from \cite{balaji},\,\S5 where the idea is 
 attributed to S. Ramanan.
Pick a curve $X_{0}$ of genus $g$ such that the Neron-Severi group
of its Jacobian, $NS(J(X_{0}))$ is $\bbZ$. By \cite{mori}
such an
$X_{0}$ exists. Pick a line bundle $L_{0}$ of degree $d$ on
$X_{0}$, and a set of $d-1$ distinct points $\chi_{0} =
\{x^1_{0},\ldots,\,x^{d-1}_{0}\}$ in $X_{0}$. One finds a
family of curves $\widetilde{X} \to T$, a line bundle $\cl$ on $\widetilde{X}$, 
and a set
of $d-1$ mutually disjoint $T$-valued points $\widetilde{\chi} =
\{\widetilde{x}^{\,1},\ldots,\,\widetilde{x}^{d-1} \}$, so that 
$(\widetilde{X},\,\cl,\,\widetilde{\chi})$ interpolates between 
$(X_{0},\,L_{0},\,\chi_{0})$ and $(X,\,L,\,\chi)$. 
To get such a triple, first observe that since the moduli space
${\cal{M}}_{g,d-1}$ of pointed curves is irreducible and
quasi-projective, we can find $(\widetilde{X}\to T,\, \widetilde{\chi})$
interpolating between $(X_0,\,\chi_0)$ and $(X,\,\chi)$.
  Let $\widetilde{\pic}^d \to T$ be the 
resulting family of degree $d$ components of the Picard groups.
Since $L_{0}$ and $L$ are points on $\widetilde{\pic}^d$, 
one can
connect them by a (possibly singular, incomplete) curve $T'$. Base change
everything to $T'$. Renaming $T'$ as $T$ and the resulting family
of pointed curves as $(\widetilde{X},\,\widetilde{\chi})$ we get a
$T$-valued point of the resulting bundle of degree $d$ components
of the Picard groups.
The line bundle $\cl$
on $\widetilde{X}$ corresponding to this section completes the triple.
We denote the specialization of
$(\widetilde{X},\,\cl,\,\widetilde{\chi})$ at $t\in T$ by $(X_t,\,L_t,\,{\chi}_t)$.
Let $t_0, t_1 \in T$ be points where 
$(X_0,\,L_0,\,{\chi}_0)$ and $(X, L, \chi)$ are realized.

The $\cs\cu_{X_t}(n,\,L_t)$ string themselves into a family $\widetilde{\cs} \to T$
(one uses Geometric Invariant Theory over the base $T$ to get $\widetilde{\cs}$.
The specializations behave well since we are working over $\bbC$). Similarly
we have a family $\widetilde{\cs^s} \to T$ specializing at $t \in T$ to
$\cs\cu^s_{X_t}(n,\,L_t)$. The intermediate Jacobians 
$J^2(\cs\cu^s_{X_t}(n,\,L_t))$ also string together into
a family of abelian varieties ${\cal{A}} \to T$. Let $\cj\to T$ be the
family $\{J(X_t)\}$ of Jacobians. Step\,1 then gives
an isomorphism of group schemes
$$
\widetilde{\vp}: \cj\lr {\cal{A}}
$$
which specializes at $t \in T$ to $\vp_{X_t,L_t,\chi_t}$. By Step\,2
we get a family of polarizations $\{\mu_t=\mu_{X_t,L_t,\chi_t}\}_{t\in T}$
on $\cj$. Since $NS(J(X_0)) = \bbZ$, therefore there exists
an integer $m \ne 0$, such that
$$
m\omega_{X_0} = \mu_{t_0}
$$
where, for any curve $C$, $\omega_C$ denotes the principal polarization
on $J(C)$. Since $\{\omega_{X_t}\}$ is a family of polarizations on
$\cj$ and since the Neron-Severi group is discrete, therefore
$$
m\omega_{X_t} = \mu_t \qquad \qquad (t \in T).
$$

Theorem\,\ref{thm:main} is now immediate.
\end{step}
\subsection{The isomorphism $\psi_{X,L,\chi}$.}\label{ss:psi}
One produces $\psi_{X,L,\chi}$ as follows\,: 
Let 
$$
\cs_1 =
\cs\cu_X(n,\,L\otimes\co_X(-D))
$$
where $D$ is the divisor
$\{x^1\}+\ldots +\{x^{d-1}\}$. Since the degree of $L\otimes\co_X(-D)$
is $1$, therefore $\cs_1$ is smooth and there exists a Poincar{\'e}
bundle $\cw$ on $X\times\cs_1$. Let $\cw_1,\ldots,\,\cw_{d-1}$ be the
$d-1$ vector bundles on $\cs_1$ obtained by restricting $\cw$ to
$\{x^1\}\times{\cs_1}=\cs_1,\ldots,\,\{x^{d-1}\}\times{\cs_1}=\cs_1$
respectively. Let $\bbP_k=\bbP(\cw_k)$, $k=1,\ldots,\,d-1$, and
$\bbP$ $(=\bbP_{X,L,\chi})$ be the product 
$\bbP_1\times_{\cs_1}\ldots\times_{\cs_1}\bbP_{d-1}$.
We will show (in \S\ref{s:hecke}) that there is a correspondence
\begin{equation}\label{eqn:hecke}
\cs_1 \stackrel{\scriptstyle{\pi}}{\lf} {\bbP} \stackrel{\scriptstyle{f}}{\lr}
\cs 
\end{equation}
where $\pi=\pi_{X,L,\chi}$ is the natural projection and $f=f_{X,L,\chi}$
is defined (via a generalized Hecke correspondence) in
\ref{ss:map-f} (see \eqref{eqn:f}).
We have isomorphisms of (integral, pure)
Hodge structures
\begin{equation}\label{eqn:tate}
H^1(X,\,\bbZ)(-1) \riso H^3(\cs_1,\,\bbZ) \riso H^3(\bbP,\,\bbZ).
\end{equation}
where the first isomorphism is that in \cite{N-R},\,p.\,392,\,Theorem\,3, and
the second is given by Leray-Hirsch.
Let $\bbP^s = f^{-1}(\cs^s)$. In \S\ref{s:hecke} (see 
Remark\,\ref{rmk:proj-bundle}, and \ref{ss:codim}) we will show

\begin{prop}\label{prop:hecke}
\begin{enumerate}
\item[(a)] If $n\ge 3$ and $g \ge 3$, the codimension of $\bbP\sm \bbP^s$
in $\bbP$ is at least $3$.
\item[(b)] The map $\bbP^s \to \cs^s$ is a 
$\bbP^{n-1}\times\ldots\times\bbP^{n-1}$
bundle, where the product is $(d-1)$-fold.
\end{enumerate}
\end{prop}

Note that if $n=2$, the codimension of $\bbP \sm \,\bbP^s$ in
$\bbP$ is $g-1$ (see \cite{balaji-thesis},\,p.\,11,\,Prop.\,7), so
that if $g\ge 4$ the codimension is at least $3$. This fact, along with
and 
Proposition\,\ref{prop:hecke} implies that the codimension of
$\bbP\sm \bbP^s$ is greater than equal to $3$ for $n,\,g$ in the
range of Theorem\,\ref{thm:main}. It then follows, from Lemma\,\ref{lem:codim}
below, that the restriction maps
\begin{align*}
H^3(\bbP,\,\bbZ) & \lr H^3(\bbP^s,\,\bbZ) \\
H^1(\bbP,\,\bbZ) & \lr H^1(\bbP^s,\,\bbZ)
\end{align*}
are isomorphisms of Hodge structures. 
Note that this means:
\begin{itemize}
\item The Hodge structure of $H^3(\bbP^s)$ is pure of 
weight $3$;
\item The cohomology group $H^1(\bbP^s,\,\bbZ) = 0$. Indeed, 
$\bbP$ is unirational (for $\cs_1$ is
--- see \cite{css-drez},\,pp.\,52--53,\,VI.B), whence $H^1(\bbP,\,\bbZ) = 0$.
\end{itemize}
We can now relate the
Hodge structures on $H^1(\cs^s)$ and $H^3(\cs^s)$ with those
on $H^1(\bbP^s)$ and $H^3(\bbP^s)$ using the map $f$ and
part (b) of Proposition\,\ref{prop:hecke}. 
For the rest of this section let $f$ also
denote the map $\bbP^s \to \cs^s$. We claim that
\begin{equation}\label{eqn:hodge}
f^*: H^3(\cs^s) \to H^3(\bbP^s)
\end{equation}
is an isomorphism of Hodge structures, modulo torsion.
This implies that
the Hodge structure on $H^3(\cs^s,\,\bbZ)$ is pure of weight $3$, a fact
that also follows from Corollary\,\ref{cor:lefschetz}.

To prove that \eqref{eqn:hodge} is an isomorphism of Hodge structures,
modulo torsion, we need:
\begin{lem}\label{lem:pi1} $\cs^s$ is simply connected.
\end{lem}
\begin{pf}
$\bbP$ is unirational, 
therefore it is simply connected \cite{pi1}. 
Since $\codim{(\bbP\sm \bbP^s)} > 1$,
it follows that $\bbP^s$ is also simply connected (purity of the
branch locus). The lemma now follows from the homotopy exact sequence
for $f$.
\end{pf}
\begin{cor} $H^1(\cs^s,\,\bbZ) = 0$.
\end{cor}
\begin{cor} $f_*\bbZ = \bbZ$, $R^1f_*\bbZ = R^3f_*\bbZ = 0$ and 
$R^2f_*\bbZ = \bbZ^{d-1}$.
\end{cor}
\begin{pf} As $\cs^s$ is simply connected, $R^if_*\bbZ$ is just the
constant sheaf associated to the $i$-th cohomology of 
$\bbP^{n-1}\times\ldots\times\bbP^{n-1}$.
\end{pf}
One can now verify \eqref{eqn:hodge} by using the Leray spectral
sequence combined with the above isomorphisms. It follows that
$H^3(\bbP^s,\,\bbZ)$ is isomorphic to the cokernel of the 
differential
$$
H^0(R^2f_*\bbZ) \to H^3(f_*\bbZ)
$$
but this vanishes mod torsion by \cite{deligne68}.
The isomorphisms \eqref{eqn:tate} and the map \eqref{eqn:hodge}, give the
desired mod-torsion isomorphism
$$
\psi_{X,L,\chi}\col H^1(X)(-1) \riso H^3(\cs^s).
$$
\begin{rem}\label{rmk:var-hodge}
This isomorphism varies well with $(X,\,L,\,\chi)$ 
as the construction
of the correspondence \eqref{eqn:hecke} will show (see 
Remark\,\ref{rmk:variation}). 
\end{rem}

Here then is the promised Lemma:
\begin{lem}\label{lem:codim} If $Y$ is a smooth projective variety,
$Z$ a codimension $k$ closed subscheme, and $U=Y\sm Z$, then
$$
H^j(Y,\,\bbZ) \riso H^j(U,\,\bbZ)
$$
for $j < 2k-1$.
\end{lem}
\begin{pf}
We have to show that $H^j_Z(X,\,\bbZ)$ vanishes for $j < 2k$. By Alexander 
duality (see for e.g. \cite{iverson},\,p.\,381,\,Theorem\,4.7) we have 
$$
H^j_Z(Y,\,\bbZ) \riso H_{2m-j}(Z,\,\bbZ),
$$
where $m=\dim{Y}$ and $H_*$ is Borel-Moore homology. Now use 
\cite{iverson},\,p.\,406,\,3.1 to conclude that the
right side vanishes if $j < 2k$ (note that $``\dim"$ in 
{\it loc.cit} is dimension as an analytic space,
and in {\it op.cit.} it is dimension as a topological (real)
manifold).
\end{pf}
\begin{rem}\label{rmk:balaji} In view of the above Lemma, it
seems that Balaji's proof of Torelli (for Seshadri's
desingularization of $\cs\cu_X(2,\,\co_X)$) does not work
for $g=3$, for in this case, the codimension of $\bbP\sm \bbP^s =2$.
(See \cite{balaji},\,top of p.\,624 and \cite{balaji-thesis},\,Remark\,9.)
\end{rem}

\subsection{The Polarization on $H^3(\cs^s)$.}\label{ss:polar}
It remains to impose a polarization on the Hodge structure of
$H^3(\cs^s)$ which varies well with $\cs^s$. Note that the map
$\psi_{X,L,\chi}$ tells us that the Hodge structure on $H^3(\cs^s)$
is pure.

One knows from the results of Drezet and Narasimhan \cite{drez-nar}, that
$\pic(\cs^s)=\bbZ$ (see p.\,89, 7.12 (especially the proof)
of {\it loc.cit.}). Moreover,
$\pic(\cs) \to \pic(\cs^s)$ is an isomorphism. Let $\xi'$ be the ample
generator of $\pic(\cs^s)$. It is easy to see that there exists a positive
integer $r$, independent of $(X,\,L)$ (with genus $X=g$), such that
$\xi = {\xi'}^r$ is very ample on $\cs$ (we are not distinguishing 
between line bundles on $\cs^s$ and their (unique) extensions to
$\cs$). Embed $\cs$ in a suitable projective space via $\xi$. Let
$e=\codim(\cs\sm \,\cs^s)$. Let $M$ be the intersection of $k = \dim\cs -e +1$
hyperplanes (in general position) with $\cs^s$. Then $M$ is smooth, projective
and contained in $\cs^s$. Let $p = \dim{\cs}$ and $H^*_c$ --- cohomology with
compact support. We then have a map
$$
l\col H^3(\cs^s) \lr H^{2p-3}_c(\cs^s)
$$
defined by
$$
x \mapsto x\cup c_1(\xi)^{p-k-3}\cup [M].
$$
If $M'$ is another $k$-fold intersection of general hyperplanes, then
$[M'] = [M]$. Hence $l$ depends only on $\xi$. According to 
Proposition\,\ref{prop:polarization} (see also Remark\,\ref{rmk:polarization}),
the pairing on $H^3(\cs^s,\,\bbC)$ given by
$$
<x,\,y> = \int_{\cs^s}l(x)\cup\,y
$$
gives a polarization on the Hodge structure of $H^3(\cs^s)$. Since
$\xi$ ``spreads" (for $\xi'$ clearly does), therefore this polarization
varies well with $\cs^s$. Then by arguments already indicated in the 
beginning of this section, this polarization is a multiple of 
principal polarization (and the integer factor is necessarily unique).
Thus one gets a natural principal polarization on $H^3(\cs^s)$.

\begin{rem}\label{rmk:saito} There is another approach to this
polarization, using Intersection Cohomology (middle perversity)
and M. Saito's theory of Hodge modules \cite{saito}. The very
ample bundle $\xi$ gives rise to Lefschetz operators 
$L^i\col IH^q(\cs) \lr IH^{q+2i}(\cs)$ (see \cite{del-beil-bern}).
Our codimension estimates (see Remark\,\ref{rmk:codim}) are such
that $IH^3(\cs) \riso H^3(\cs^s)$ and $IH^1(\cs) = H^1(\cs^s) = 0$.
The group $IH^3(\cs)$ has a pairing on it given by
$$
<\alpha,\,\beta> = \int_SL^{p-3}\alpha\cup\beta
$$
where $\int_S(\_)\cup\beta\col IH^{2p-3}(\cs)\to \bbC$ is the map
given by the Poincar{\'e} duality pairing between $IH^{2p-3}(\cs)$
and $IH^3(\cs)$. According to M. Saito \cite{saito},\,5.3.2, this
gives a polarization on the Hodge structure of $IH^3(\cs)$ (since
all classes in $IH^3(\cs)$ are primitive). This polarization 
translates to one on $H^3(\cs^s)$. A little thought shows (say by
desingularizing $\cs$) that the pairing on $H^3(\cs^s)$ is
$$
<x,\,y> = \int_{\cs^s}c_1(\xi)^{p-3}\,\wedge\,x\,\wedge\,y.
$$
Here, on the right side, we are using De Rham theory, and replacing the
various elements in cohomology by forms which represent them. The
integral above is the usual integral of forms. Note that we
could not have defined the pairing by the above formula, for we
have no {\it a priori} guarantee that the right side (which is
an integral over an open manifold) is finite.
\end{rem}

\section{The correspondence variety $\bbP$}\label{s:hecke}
In this section we define the map $f\col \bbP \to \cs$ and prove 
Proposition \ref{prop:hecke}.
\subsection{The map $f\col \bbP \to \cs$.}\label{ss:map-f}
 We need some notations\,:
\begin{itemize}
\item For $1\le k\le d-1$, $\pi_k\col \bbP \to \bbP_k$ is 
the natural projection;
\item $\im\col Z \hookrightarrow X$ is the reduced subscheme defined
by $\chi=\{x^1,\ldots,\,x^{d-1}\}$.
\item $\im_k\col Z_k \hookrightarrow X$, the reduced scheme defined by
$\{x_k\}$, $k= 1,\ldots,\,d-1$.
\item For any scheme $S$,
\begin{enumerate}
\item[(i)] $p_S\col X\times\,S \to S$ and $q_S\col X\times\,S \to X$ are
the natural projections;
\item[(ii)] $Z^S = q_S^{-1}(Z)$;
\item[(iii)] $Z_k^S = q_S^{-1}(Z_k)$, $k=1,\ldots,\,d-1$.
Note that $Z^S_k$ can be identified canonically with $S$.
\end{enumerate}
\end{itemize}

We will show --- in \ref{ss:u-exact} --- that there is an exact sequence
\begin{equation}\label{eqn:u-exact}
0 \lr (1\times\pi)^*\cw \lr \cv \lr \ct_0 \lr 0
\end{equation}
on $X\times\,\bbP$, with $\cv$ a vector bundle on $X\times\bbP$
and $\ct_0$ a line bundle {\it on the subscheme $Z^\bbP$}, 
which is universal
in the following sense\,: If $\psi\col S \to \cs_1$ is a $\cs_1$-scheme
and we have an exact sequence
\begin{equation}\label{eqn:s-exact}
0 \lr (1\times\psi)^*\cw \lr \ce \lr \ct \lr 0
\end{equation}
on $X\times\,S$, with $\ce$ a vector bundle on $X\times\,S$
and $\ct$ a line bundle {\it on the subscheme} $Z^S$, then
there is a unique map of $\cs_1$-schemes
$$
g\col S \lr \bbP
$$
such that, 
$$
(1\times\,g)^*\eqref{eqn:u-exact} \equiv  \eqref{eqn:s-exact}.
$$
The $\equiv$ sign above means that the two exact sequences
are isomorphic, and the left most isomorphism 
$(1\times{g})^*\comp(1\times{\pi})^* \riso (1\times\psi)^*$ is the canonical
one.

There is a way of interpreting this universal property in terms of
quasi-parabolic bundles (see \cite{mehta-css},\,p.\,211--212,\,Definition\,1.5,
for the definitions of quasi-parabolic and parabolic bundles). Taking
$\chi$ as our collection of parabolic vertices, we can introduce a
quasi-parabolic datum on $X$ by attaching the flag type $(1,\,n-1)$ to
each point of $\chi$. From now onwards {\it quasi-parabolic structures
will be with respect to this datum and on vector bundles of rank $n$ and
determinant $L$}. 
One observes that for a vector bundle $V$ (of rank $n$ and determinant
$L$), a surjective map $V\twr \co_Z$ determines a unique quasi-parabolic
structure, and two such surjections give the same quasi-parabolic strcuture
if and only if they differ by a scalar multiple.
The above mentioned universal property says that $\bbP$
is a (fine) moduli space for quasi-parabolic bundles. More precisely, the
family of quasi-parabolic structures
$$
\cv \twr \ct_0
$$
parameterized by $\bbP$ is universal for families of quasi-parabolic bundles
$$
\ce \twr \ct
$$
parameterized by $S$, whose kernel is a family of semi-stable bundles. The
points of $\bbP$ parameterize quasi-parabolic structures $V \twr \co_Z$
whose kernel is semi-stable.

Let $\balp = (\alpha_1,\,\alpha_2)$, where $0 < \alpha_1 < \alpha_2 <1$,
and let $\D = \D_{\balp}$ be the parabolic datum which attaches to each
parabolic vertex (of our quasi-parabolic datum) weights $\alpha_1, \alpha_2$.
We can choose $\alpha_1$ and $\alpha_2$ so small that
\begin{itemize}
\item a parabolic semi-stable bundle is parabolic stable\,;
\item if $V$ is stable, then every parabolic structure on $V$ is parabolic
stable\,;
\item the underlying vector bundle of a parabolic stable bundle is
semi-stable in the usual sense\,;
\item if $V \twr \co_Z$ is parabolic stable, then the kernel $W$ is
semi-stable.
\end{itemize}
Showing the above involves some very elementary calculations. Denote
the resulting moduli space of parabolic stable bundles 
$\cs\cu_X(n,\,L,\,\D)$.

Let $\bbP^{ss}\subset \bbP$ be the locus on which $\cv$ consists
of parabolic semi-stable (=parabolic stable) bundles. One checks that
$\bbP^{ss}$ is an open subscheme of $\bbP$ (this involves two things\,:
(i) knowing that the scheme $\widetilde{R}$ of \cite{mehta-css},\,p.\,226
has a local universal property for parabolic bundles and (ii) knowing that
the scheme ${\widetilde{R}}^{ss}$ of {\it loc.cit.} is open).

Clearly $\bbP^{ss}$ is non-empty --- in fact if $V$ is stable of rank $n$
and determinant $L$, then any parabolic structure on $V$ is parabolic
stable (see above). We claim that $\bbP^{ss} \simeq \cs\cu_X(n,\,L,\,\D)$.
To that end, let $S$ be a scheme, and
\begin{equation}\label{eqn:par-fly}
\ce \twr \ct
\end{equation}
a family of parabolic stable bundles parameterized by $S$. The kernel
$\cw'$ of \eqref{eqn:par-fly} is a family of stable bundles of rank $n$
and determinant $L\otimes\co_X(-D)$. Since $\cs_1$ is a fine moduli space,
we have a unique map $g\col \cs\to \cs_1$ and a line bundle $\xi$ on $S$
such that $(1\times{g})^*\cw = \cw'\otimes\,p_S^*\xi$. By doctoring
\eqref{eqn:par-fly} we may assume that $\xi=\co_S$. The universal
property of the exact sequence \eqref{eqn:u-exact} on $\bbP$ then
gives us a unique map
$$
g\col\,S \lr \bbP
$$
such that $(1\times{g})^*\eqref{eqn:u-exact}$ is equivalent to
$$
0 \lr \cw' \lr \ce \lr \ct \lr 0.
$$
Clearly $g$ factors through $\bbP^{ss}$. This proves that
$\bbP^{ss}$ is $\cs\cu_X(n,\,L,\,\D)$. However, $\cs\cu_X(n,\,L,\,\D)$
is a projective variety (see \cite{mehta-css},\,pp.\,225--226,\,Theorem\,4.1),
whence we have 
$$
\bbP = \cs\cu_X(n,\,L,\,\D).
$$
It follows that $\cv$ consists of parabolic stable bundles, and hence
of (usual) semi-stable bundles (by our choice of $\balp$). Since $\cs$
is a coarse moduli space, we get the map
\begin{equation}\label{eqn:f}
f\col\,\bbP \lr \cs .
\end{equation}

\begin{rem}\label{rmk:hecke}
Note that the parabolic structure $\D$ is something of a red herring.
In fact $\cs\cu_X(n,\,L,\,\D)$ parameterizes quasi-parabolic structures
$V \twr \co_Z$, whose kernel is semi-stable
(cf. \cite{mehta-css},\,p.\,238,\,Remark\,(5.4), where this point is made
for $n=2, d=2$). The space $\bbP$ should be thought of as the correspondence
variety for a certain Hecke correspondence (cf. \cite{N-R-hecke}).
\end{rem}

\begin{rem}\label{rmk:proj-bundle} 
Let $V$ be a stable bundle of rank $n$, with $\det{V}=L$, so that
(the isomorphism class of) $V$ lies in $\cs^s$. Since any parabolic
structure on $V$ is parabolic stable (by our choice of $\balp$), therefore
we see that $f^{-1}(V)$ is canonically isomorphic to
$\bbP(V_{x^1}^*)\times\ldots\times\bbP(V_{x^{d-1}}^*)$. \footnote{One can
be more rigorous. Identifying $Z_k^{\bbP}$ with ${\bbP}$ for each
$k = 1,\ldots,\,d-1$, we see that restricting the universal exact
sequence to $Z_k^{\bbP}$ gives us $d-1$ quotients 
$\co_{\bbP}\otimes_{\bbC}V_{x^k} \twr \ct_0|Z^{\bbP}_k$. Let $S$
be a scheme which has $d-1$ quotients $\co_S\otimes_{\bbC}V_{x^k}\twr \cl_k$
$k = 1,\ldots,\,d-1$, on it, where the $\cl_k$ are line bundles. These
quotients extend to a family of parabolic structures
$q_S^*V \twr \ct$ (on $V$) parameterized by $S$ in a unique way. The universal
property of the exact sequence \eqref{eqn:u-exact} gives us a map
$S \to \bbP$, and this map factors through $f^{-1}(V)$.} This gives us
part (b) of Proposition\,\ref{prop:hecke}, for it is not hard to see
that $\bbP^s \to \cs^s$ is smooth (examine the effect on the tangent space
of each point on $\bbP^s$).
\end{rem}
\subsection{Codimension estimates.}\label{ss:codim}  We wish to estimate 
$\codim{(\bbP\sm \bbP^s})$. For any vector bundle $E$ on $X$, let
$\mu(E)=\rank{E}/\deg{E}$. Let $\mu=d/n$. Let $V \twr \co_Z$ be
a parabolic bundle in $\bbP\sm \bbP^s$. Then we have a filtration
(see \cite{css-drez},\,p.\,18,\,Th{\'e}or{\`e}me\,10)
$$
0=V_{p+1} \subset V_p \subset \ldots \subset V_0=V
$$
such that for $0\le i \le p$, $G_i=V_i/V_{i+1}$ is stable and $\mu(G_i)=\mu$.
Moreover (the isomorphism class of) the vector bundle $\bigoplus{G_i}$ depends
only upon $V$ and not on the given filtration. We wish to count the
number of moduli at $[V\overset{\theta}{\twr} \co_Z] \in \bbP\sm \bbP^s$.
There are three sources\,:
\begin{enumerate}
\item[a)] The choice of $\bigoplus_{i=0}^pG_i$\,;
\item[b)] Extension data\,;
\item[c)] The choice of parabolic structure $V\overset{\theta}{\twr}\co_Z$,
for fixed semi-stable $V$.
\end{enumerate}
The source c) is the easiest to calculate --- there is a codimension
one subspace at each parabolic vertex, contributing
$$
(n-1)(d-1)
$$
moduli.

Let $n_i=\rank{G_i}$.

The number of moduli arising from a) is evidently
$$
\sum_{i=0}^p(n_i^2-1)(g-1) + pg .
$$
Indeed, the bundles $G_i$ have degree $n_i\mu$ and the product of their
determinants must be $L$. They are otherwise unconstrained.

It remains to estimate the number of moduli arising from extension
data. Each extension
$$
0 \lr V_{i+1} \lr V_i \lr G_i \lr 0 \qquad i=0,\ldots,\,p
$$
determines a class in $H^1(X,\,G_i^*\otimes{V_{i+1}})$. Note that
\begin{equation*}
\begin{split}
h^0(G_i^*\otimes{V_{i+1}}) & = \dim{\Hom_{\co_X}(G_i,\,V_{i+1})} \\
 & \le \sum_{j>i}\Hom_{\co_X}(G_i,\,G_j) \\
 & \le p-i
\end{split}
\end{equation*}
by the sub-additivity of $\dim{\Hom(G_i,\,\_)}$ and the stability of $G_i$.
By the Riemann-Roch theorem 
\begin{equation*}  
\begin{split}
h^1(G_i^*\otimes{V_{i+1}}) & = h^0(G_i^*\otimes{V_{i+1}}) - 
n_i(n_{i+1}+\ldots +n_p)(1-g) \\
 & \le (p-i) - n_i(n_{i+1}+\ldots n_p)(1-g).
\end{split}
\end{equation*}  
The isomorphism class of $V_i$ depends only on a scalar multiple of the
extension class. Therefore the number of moduli contributed by extensions
is
\begin{equation*}
\begin{split}
\sum_{i=0}^p\left [h^1(G_i^*\otimes{V_{i+1}} -1 \right]  & \le
\sum_{i=0}^p\left[p-i-n_i(n_{i+1}+\ldots n_p)(1-g) \right] - (p+1) \\
 & = \dfrac{p(p+1)}{2} - 
\sum_{i=0}^{p-1}n_i(n_{i+1}+\ldots +n_p)(1-g) - (p+1) \\
 & = \dfrac{(p+1)(p-2)}{2} - \sum_{i<j}n_in_j(1-g).
\end{split}
\end{equation*}
Adding the contributions from a), b) and c) and subtracting from
$$
\dim\bbP=(n^2-1)(g-1)+(n-1)(d-1)$$ 
we get
\begin{equation*}
\begin{split}
\codim(\bbP\sm \bbP^s) & \ge (n^2-1)(g-1) - \sum_{i=o}^p(n_i^2-1)(g-1) - pg \\
 & \quad - \sum_{i<j}n_in_j(g-1) - \dfrac{(p+1)(p-2)}{2} \\
& = \sum_{i<j}n_in_j(g-1) - \dfrac{(p-1)(p+2)}{2} \\
& = B \qquad (\text{say}).
\end{split}
\end{equation*}
Now,
$\sum_{i<j}n_in_j \ge {p(p+1)}/{2}$, therefore
$$
B \ge \dfrac{p(p+1)}{2}(g-1) - \dfrac{(p+2)(p-1)}{2}.
$$
It follows that $B\ge 3$ whenever $p\ge 2$ {\it and} $g\ge 3$. If $p=1$ and
$n \ge 3$, then 
$$
B/(g-1) = \sum_{i<j}n_in_j \ge 2
$$
and one checks that $B\ge 3$ whenever $g\ge 3$. 
Proposition\,\ref{prop:hecke}(a) may now be considered as proved. 

\begin{rem}\label{rmk:codim} We could use similar techniques
to estimate $\codim{(\cs\sm \cs^s)}$, but our task is made
easier by the exact answers in \cite{css-drez},\,p.\,48,\,A. For just
this remark, assume $d > n(2g-1)$, and let $a = (n,\,d)$. Then
$a\ge 2$. Let $n_0 = n/a$. Then according to {\it loc.cit.},
\begin{equation*}
\codim{(\cs\sm \cs^s)}=
\begin{cases}
(n^2-1)(g-1) -
\dfrac{n^2}{2}(g-1) -2 +g & \text{if $a$ is even} \\
 & {} \\
(n^2-1)(g-1) - \dfrac{n^2 + n_0^2}{2}(g-1) -2 + g & \text{if $a$ is odd}.
\end{cases}
\end{equation*}
It now follows that
$$
\codim{(\cs\sm \cs^s)} > 5
$$
if $n,\,g$ are in the range of Theorem\,\ref{thm:main}.
\end{rem}
\subsection{The universal exact sequence on $X\times\bbP$.}\label{ss:u-exact}
We begin by reminding the reader of some elementary facts from commutative
algebra. If $A$ is a ring (commutative, with $1$), $t\in A$ a non-zero
divisor, and $M$ an $A$-module, then each element $m_0\in M$ gives
rise to an equivalence class of extensions
\begin{equation}\label{eqn:exact-mod}
0 \lr M \lr E_{m_0} \lr A/tA \lr 0
\end{equation}
where $E_{m_0} = \left(A\bigoplus{M}\right)/A(t,\,m_0)$, and the
arrows are the obvious ones. Moreover, if $m_0 - m_1 \in tM$, say
$$
m_0 - m_1 = tm'
$$
then the extension given by $m_0$ is equivalent to that given by $m_1$.
In fact, one checks that
\begin{equation}\label{eqn:patching}
\begin{split}
E_{m_0} & \riso E_{m_1} \\
(a,\,m) & \mapsto (a,\,m-am')
\end{split}
\end{equation}
gives the desired equivalence of extensions. This is another way of
expressing the well known fact that each element of
$M/tM = \ext^1(A/t,\,M)$ gives rise to an extension. 

One globalizes to get the following\,: Let $S$ be a scheme, 
$T\overset{\im}{\hookrightarrow} S$ a closed immersion, $\cf$ a
quasi-coherent $\co_S$-module, $U$ an open neighbourhood of $T$ in $S$,
and $t\in \G(U,\,\co_S)$ an element which defines $T \hookrightarrow U$,
and which is a non-zero divisor for $\G(V,\,\co_S)$
for any open $V \subset U$. Then every global section $s$ of 
$\im^*\cf = \cf\otimes\co_T$ gives rise to an equivalence class of
extensions
\begin{equation}\label{eqn:patched}
0 \lr \cf \lr \ce \lr \co_T \lr 0.
\end{equation}
Indeed, we are reduced immediately to the case $S = U$. We build up
exact sequences \eqref{eqn:exact-mod} on each affine open subset 
$W \subset S$, by picking a lift $\tilde{s_W}\in\G(W,\,\cf)$ of
$s\,|\,W$. One patches together these exact sequences via 
\eqref{eqn:patching}.

Now consider $\bbP = \bbP_1\times_{\cs_1}\ldots\times_{\cs_1}\bbP_{d-1}$. 
For each
$k = 1,\ldots,\,d-1$, let $p_k\col\,\bbP_k \to \cs_1$ be the natural
projection. We have a universal exact sequence
$$
0 \lr \co(-1) \lr p_k^*\cw_k \lr B \lr 0
$$
whence a global section $s_k\in\G(\bbP_k,\,p_k^*\cw_k(1))$. 
However, note that
$$
p_k^*\cw_k = (1\times p_k)^*\cw\,|\,Z_k^{\bbP_k}
$$
where we are identifying $Z_k^{\bbP_k}$ with $\bbP_k$.
By \eqref{eqn:patched} we get exact sequences
$$
0 \lr (1\times\pi)^*\cw \lr \cv_k \lr
\co_{Z_k^\bbP}\otimes{L_k} \lr 0 
$$
where $L_k$ is the line bundle obtained by pulling up $\co_{\bbP_k}(-1)$.
It is not hard to see
that $\cv_k$ is a family of vector bundles
parameterized by $\bbP$. Glueing these sequences together --- the $k$-th
and the $l$-th agree outside $Z^{\bbP}_k$ and $Z^{\bbP}_l$ --- we
obtain \eqref{eqn:u-exact}.

Now suppose we have a $\cs_1$-scheme $\psi\col\,S \to \cs_1$ and the
exact sequence \eqref{eqn:s-exact}
$$
0 \lr (1\times\psi)^*\cw \lr \ce \lr \ct \lr 0
$$
on $X\times{S}$. Restricting \eqref{eqn:s-exact} to $Z^S_k$ ($1\le k \le d-1$)
one checks that the kernel of $(1\times\psi)^*\cw\,|\,Z^S_k \to \ce\,|Z_k^S$
is a line bundle $\cl_k$. Identifying $Z^S_k$ with $S$, we see that
$(1\times\,\psi)^*\cw\,|\,Z^S_k = \psi^*\cw_k$. Thus $\cl_k$ is a line
sub-bundle of $\psi^*\cw_k$. By the universal property of $\bbP_k$, we
see that we have a unique map of $\cs_1$-schemes
$$
g_k\col\,S \lr \bbP_k
$$
such that $\co(-1)$ gets pulled back to $\cl_k$. The various $g_k$
give us a map
$$
g\col S \lr \bbP
$$
One checks that $g$ has the required universal property. The uniqueness
of $g$ follows from the uniqueness of each $g_k$.

\begin{rem}\label{rmk:variation} It is clear from the construction
that the map
$$
f = f_{X,L,\chi}\col\,\bbP_{X,L,\chi} \lr \cs\cu_X(n,\,L)
$$
varies well with $(X,\,L,\,\chi)$. This implies that the correspondence
\eqref{eqn:hecke} also varies well with $(X,\,L,\,\chi)$ and hence so
do $\psi_{X,L,\chi}$ and $\vp_{X,L,\chi}$.
\end{rem}

\section{Polarizations}\label{s:polarization}
Let $Y$ be an $m$-dimensional projective variety.
Suppose that $U$ is a smooth Zariski open subset. 
One then has the following version of the
Lefschetz theorem.

\begin{thm}\label{thm:lefschetz} If $H$ is a hyperplane section
of $Y$ such that $U\cap H$ is non-empty, then
$$
H^i(U,\,\bbQ) \to H^i(U\cap{H},\,\bbQ)
$$
is an isomorphism for $i < m-1$ and injective when $i = m-1$.
\end{thm}
\begin{pf} We need some results involving Verdier duality. The
standard references are \cite{borel} and \cite{iverson}. Let
$S$ be an analytic space and $p_S$ the map from $S$ to a point.
For $\cf\in D^b_{const}(S,\,\bbQ)$ (the derived category of
bounded complexes of $\bbQ_S$-sheaves whose cohomology sheaves
are $\bbQ_S$-constructible), set
$$
D_S(\cf) = \bbR\ch\morph_S(\cf,\,p_S^!\bbQ). 
$$
We then have by Verdier duality
\begin{equation}\label{eqn:verdier}
\bbH^i(S,\,\cf) \overset{\sim}{\lr} \bbH^{-i}(S,\,D_S(\cf))^*.
\end{equation}
Here $\bbH^*$ denotes ``hypercohomology".

For an open immersion $h\col S' \hookrightarrow S$, one has canonical
isomorphisms
\begin{align}\label{eqn:*!}
\bbR{h}_*D_{S'}\cg & \overset{\sim}{\lr} D_S(h_!\cg) \\
\intertext{and}
\label{eqn:!*}
{\bbR}h_!D_{S'}\cg & \overset{\sim}{\lr} D_S(\bbR{h_*}\cg)
\end{align}
Here $\cg\in D^b_{const}(S',\,\bbQ)$.
The first isomorphism is easy (using Verdier duality for the
map $h$) and the second follows from the first and from the
fact that $D_{S'}$ is an involution. We have used (throughout)
the fact that $h_!$ is an exact functor.

If $S$ is smooth, we have
\begin{equation}\label{eqn:dim}
p_S^!\bbQ = \bbQ_S[2\dim{S}].
\end{equation}
In order to prove the theorem, let $V=U\sm H$ and $W=Y\sm H$. We then
have a cartesian square
$$
\begin{array}{ccc}
V & \stackrel{\scriptstyle{{\im}'}}{\lr} & U \\
\vcenter{%
\llap{$\scriptstyle{{\jm}'}$}}\Big\downarrow & 
& \Big\downarrow\vcenter{%
\rlap{$\scriptstyle{\jm}$}} \\
W & \underset{\scriptstyle{\im}}{\lr} & Y
\end{array}
$$
where each arrow is the obvious open immersion.
We have, by \eqref{eqn:*!} and \eqref{eqn:!*}, the identity
\begin{equation}\label{eqn:Dji}
{\jm}_!\bbR{\im}'_*D_V\bbQ_V = D_Y(\bbR\jm_*{\im}'_!\bbQ_V).
\end{equation}
Consider the exact sequence of sheaves
$$
0 \lr {\im}'_!\bbQ_V \lr \bbQ_U \lr g_*\bbQ_{H\cap U} \lr 0
$$
where $g\col H\cap{U}\to U$ is the natural closed immersion. It suffices
to prove that $H^i(U,\,{\im}'_!\bbQ_V)=0$ for $i\le m-1$. Now,
$$
H^i(U,\,{\im}'_!\bbQ_V)=\bbH^i(Y,\,\bbR\jm_*{\im}'_!\bbQ_V).
$$
Using \eqref{eqn:verdier}, \eqref{eqn:Dji} and \eqref{eqn:dim}, the
above is dual to
\begin{align*}
\bbH^{-i}(Y,\,\jm_!\bbR{\im}'_*D_V\bbQ_V) 
 & = \bbH^{2m-i}(Y,\,\jm_!\bbR{\im}'_*\bbQ_V) \\
\intertext{But $\jm_!\bbR{\im}'_* = \bbR\im_*{\jm}'_!$, and hence the above is}
 & = \bbH^{2m-i}(Y,\,\bbR{\im}_*({\jm}'_!\bbQ_V)) \\
 & = \bbH^{2m-i}(W,\,{\jm}'_!\bbQ_V) \\
 & = H^{2m-i}(W,\,{\jm}'_!\bbQ_V).
\end{align*}
Now, $W$ is an affine variety, and therefore, according to
M. Artin, its constructible cohomological
dimension is less than or equal to its dimension \cite{artin}.
Consequently, the above chain of equalities vanish whenever
$i<m$ (see also \cite{gor-mac}).
\end{pf}

We immediately have:
\begin{cor}\label{cor:lefschetz} Let $e=\codim(Y\sm U)$.
For $i< e-1$, the Hodge structure 
$H^i(U)$ is pure of weight $i$.
\end{cor}
\begin{pf} This is true if $U$ is projective. In general proceed 
using Bertini's theorem, induction, Theorem\,\ref{thm:lefschetz} and
the fact that submixed Hodge structures of pure Hodge structures are
pure \cite{deligne-hodge}.
\end{pf}

Let $i\in\bbN$ and $\cl$ a line bundle on $Y$ be such that
\begin{enumerate}
\item[(a)] $H^j(U,\,\bbQ)=0$ for $j = i-2,\,i-4,\ldots $\,;
\item[(b)] $i < e-1$\,;
\item[(c)] $\cl$ is very ample.
\end{enumerate}
\begin{rem}\label{rmk:polarization} Note that if $Y=\cs$, $U=\cs^s$, 
then $i=3$ and $\cl=\xi$ ($\xi$= the very ample bundle
of \ref{ss:polar}) satisfy the above conditions
by the results of \ref{ss:psi} and Remark\,\ref{rmk:codim}. 
\end{rem}

Let $M$ be the intersection of $k=m-e+1$  hyperplanes in general position.
Then $M$ is a smooth variety contained in $U$. Let
$$
l\col H^i(U) \lr H^{2m-i}_c(U)
$$
be the composite of 
\begin{equation*}
\begin{split}
H^i(U) & \lr H^i(M) \\
 & \lr H^{2m-2k-i}(M) \\
 & \lr H^{2m-i}_c(U)
\end{split}
\end{equation*}
where the first map is restriction, the second is 
``cupping with $c_1(\cl)^{m-k-i}$" and the third is the Poincar{\'e} dual
to restriction. The map $l$ is also described as
$$
x \mapsto x\cup c_1(\cl)^{m-k-i}\cup [M].
$$
One then has (easily)
\begin{lem} If $M'$ is another $k$-fold intersection of general hyperplanes,
then $[M'] = [M]$. Therefore $l$ depends only on $\cl$.
\end{lem}
\begin{prop}\label{prop:polarization} The pairing
$$
<x,\,y> = \int_Ul(x)\cup y
$$
on $H^i(U,\,\bbC)$
gives a polarization on the pure Hodge structure $H^i(U)$.
\end{prop}
\begin{pf} By Theorem\,\ref{thm:lefschetz}, we have an isomorphism
$$
r\col H^i(U) \lr H^i(M).
$$
The latter Hodge structure carries a polarization given by
$$
<\alpha,\,\beta> = \int_Mc_1(\cl)^{m-k-i}\cup\alpha\cup\beta.
$$
The conditions on $i$ and the Hodge-Riemann bilinear relations
on the primitive part of $H^i(M,\,\bbC)$, assure us that the above is indeed
a polarization (see \cite{grif-period} or Chap.\,V,\,\S6
of \cite{wells}). In fact, our conditions on $i$ imply that the
primitive part of $H^i(M)$ is everything.
This translates to a polarization on $H^i(U)$ given by
$$
<x,\,y> = \int_Ul(x)\cup y.
$$
This gives the result.
\end{pf}
\begin{ack} We wish to thank Prof.\,M.\,S. Narasimhan and Prof.\,C.\,S.
Seshadri for their encouragement and their help. Thanks to V. Balaji,
L. Lempert, N. Raghavendra and P.\,A.Vishwanath for helpful discussions.
Balaji made us aware of the problem, and generously discussed his proof
(in \cite{balaji}) of the Torelli theorem for Seshadri's desingularization 
of $\cs\cu_X(2,\,\co_X)$. The second author gratefully acknowledges
the four wonderful years he spent at the SPIC Science Foundation, Madras. 
\end{ack}

%\bibliography{torelli}

\bibliographystyle{plain}
\end{document}